# Hankel determinants of some polynomial sequences

*Johann Cigler*

Fakultät für Mathematik, Universität Wien

*johann.cigler@univie.ac.at*

**Abstract**

We give simple new proofs of some Catalan – Hankel determinant evaluations by Ömer Eğecioğlu and Aleksandar Cvetković, Predrag Rajković and Miloš Ivković and prove analogous results for sums of moments of "symmetric" orthogonal polynomials.

**Introduction**

Let $(a(n))_{n \geq 0}$ be a given sequence and let $r(n,x) = \sum_{k=0}^{n} a(n-k)x^k$ or $r(n,x) = xa(n) - a(n+1)$. Consider the Hankel matrices $A_n = \left(a(i+j)\right)_{i,j=0}^{n}$, $R_n = \left(r(i+j,x)\right)_{i,j=0}^{n}$ and let $d(n,x) = \det R_n$. Suppose that $d(n) = \det A_n \neq 0$ for all $n$. We want some information about the ratio $\dfrac{d(n,x)}{d(n)}$.

Ömer Eğecioğlu [10] has shown that if we choose for $a(n)$ the Catalan numbers $C_n = \dfrac{1}{n+1}\binom{2n}{n}$ or the central binomial coefficiens $B_n = \binom{2n}{n}$ and let $r(n,x) = \sum_{k=0}^{n} a(n-k)x^k$, then the quotient of the determinants $\dfrac{d(n,x)}{d(n)} = \dfrac{d(n,x)}{d(n,0)}$ is given by

$$\frac{d(n,x)}{d(n,0)} = \sum_{i=0}^{n} (-1)^i \binom{n+i}{n-i} x^i.$$

Aleksandar Cvetković, Predrag Rajković and Miloš Ivković [9] have shown that $\det\left(C_{i+j} + C_{i+j+1}\right)$ is a Fibonacci number.

We show more generally that analogous results hold if the numbers $a(n)$ are moments of "symmetric" orthogonal polynomials $p(n,x)$. By symmetric we mean that they satisfy a recurrence of the form $p_n(x) = xp_{n-1}(x) - t_{n-2}p_{n-2}(x)$. Some of these results have also been obtained with other methods in [3], [11] and [12].

# 1. The polynomials $r(n,x) = \sum_{k=0}^{n} a(n-k)x^k$.

We start with the obvious fact that $\dfrac{d(n,x)}{d(n,0)} = \det\left(R_n A_n^{-1}\right)$.

A short inspection shows that

$R_n A_n^{-1} = \left(h_n(i,j)\right)_{i,j=0}^{n}$ with $h_n(i,0) = x^i h_n(0,0)$ and $h_n(i,j) = x^i h_n(0,j)$ for $1 \leq i < j$ and $h_n(i,j) = x^{i-j} + x^i h_n(0,j)$ if $i \geq j \geq 1$.

For let $A_n^{-1} = \left(u(i,j)\right)_{i,j=0}^{n}$. Then $\sum_{k=0}^{n} a(i+k)u(k,j) = [i=j]$ and $h_n(i,j) = \sum_{k=0}^{n} r(i+k,x)u(k,j)$.

Note that $h_n(0,j) = \sum_{k=0}^{n} r(k,x)u(k,j)$. Therefore

$$h_n(i,0) = \sum_{k=0}^{n} r(i+k,x)u(k,0) = \sum_{k=0}^{n} u(k,0)\left( x^i r(k,x) + \sum_{\ell=0}^{i-1} a(k+i-\ell)x^\ell \right) = x^i h_n(0,0).$$

If $1 \leq i < j$ we get in the same way

$$h_n(i,j) = \sum_{k=0}^{n} r(i+k,x)u(k,j) = \sum_{k=0}^{n} u(k,j)\left( x^i r(k,x) + \sum_{\ell=0}^{i-1} a(k+i-\ell)x^\ell \right) = x^i h_n(0,j).$$

For $i \geq j \geq 1$ the second sum for $\ell = i-j$ adds $x^{i-j}$.

For example for $n=3$ the matrix is (we write $h(j)$ in place of $h_3(0,j)$)

$$R_3 A_3^{-1} = \begin{pmatrix} h(0) & h(1) & h(2) & h(3) \\ xh(0) & 1+xh(1) & xh(2) & xh(3) \\ x^2 h(0) & x + x^2 h(1) & 1 + x^2 h(2) & x^2 h(3) \\ x^3 h(0) & x^2 + x^3 h(1) & x + x^3 h(2) & 1 + x^3 h(3) \end{pmatrix}.$$

The special form of the matrix implies that

$$\det\left(R_n A_n^{-1}\right) = h_n(0,0). \tag{1}$$

Consider polynomials $p_n(x)$ which satisfy a recurrence of the form

$$p_n(x) = xp_{n-1}(x) - t_{n-2} p_{n-2}(x) \tag{2}$$

for some numbers $t_n \neq 0$ with initial values $p_{-1}(x) = 0$ and $p_0(x) = 1$. These polynomials have the form

$$p_n(x) = \sum_{k=0}^{\lfloor \frac{n}{2} \rfloor} (-1)^k v(n,k) x^{n-2k} \tag{3}$$

and are orthogonal with respect to some linear functional $\Lambda$, i.e. $\Lambda(p_n p_m) = 0$ for $n \neq m$. We call them symmetric orthogonal polynomials. This linear functional is uniquely determined by $\Lambda(p_n) = [n = 0]$. If the sequence $t = (t_n)$ is given we say that the polynomials $p_n(x)$ are associated with $t$.

Throughout this paper we assume that $a(n)$ is of the form

$$a(n) = \Lambda(x^{2n}) \tag{4}$$

for a sequence of symmetric orthogonal polynomials. We shortly call $a(n)$ symmetric moments. For more information on orthogonal polynomials we refer to [14].

Continuing the recurrence
$$p_{2n+1}(x) = xp_{2n}(x) - t_{2n-1}p_{2n-1}(x) = xp_{2n}(x) - t_{2n-1}(xp_{2n-2}(x) - t_{2n-3}p_{2n-3}(x)) = \cdots$$
we see that
$$\frac{p_{2n+1}(x)}{x} = p_{2n}(x) - t_{2n-1}p_{2n-2}(x) + t_{2n-1}t_{2n-3}p_{2n-4}(x) + -\cdots + (-1)^n t_{2n-1}t_{2n-3}\cdots t_1 p_0(x).$$
This implies that

$$\Lambda\left(\frac{p_{2n+1}(x)}{x}\right) = (-1)^n t_{2n-1}t_{2n-3}\cdots t_1.$$

By orthogonality

$$\Lambda\left(x^{2k-1}p_{2n+1}(x)\right) = 0$$

for $1 \leq k \leq n$.
Therefore
$$M_n(x) = \frac{(-1)^n}{t_{2n-1}t_{2n-3}\cdots t_1} \frac{p_{2n+1}(x)}{x} = \frac{(-1)^n}{t_{2n-1}t_{2n-3}\cdots t_1} \sum_{j=0}^{n} (-1)^k v(2n+1,k) x^{2n-2k} \tag{5}$$

satisfies

$$\Lambda\left(x^{2m} M_n(x)\right) = [m = 0] \tag{6}$$

for $0 \leq m \leq n$.

This means that

$$\sum_{k=0}^{n} (-1)^{n+k} \frac{1}{t_{2n-1}t_{2n-3}\cdots t_1} v(2n+1,k) a(n+m-k) = [m=0]. \tag{7}$$

Thus the transpose of the first column of $A_n^{-1}$ is

$$\left(u(0,0), u(1,0), \cdots, u(n,0)\right) = \frac{1}{t_{2n-1}t_{2n-3}\cdots t_1}\left((-1)^j v(2n+1, n-j)\right) \text{ and (1) implies that}$$

$$\det\left(R_n A_n^{-1}\right) = h_n(0,0) = \frac{(-1)^n}{t_{2n-1}t_{2n-3}\cdots t_1}\sum_{k=0}^{n}(-1)^k v(2n+1,k) r(n-k,x).$$

We state these results in the following

**Lemma**

Let $p_n(x) = \sum_{k=0}^{n}(-1)^k v(n,k) x^{n-2k}$ satisfy the recurrence $p_n(x) = x p_{n-1}(x) - t_{n-2} p_{n-2}(x)$ with $p_{-1}(x) = 0$ and $p_0(x) = 1$ and let $a(n) = \Lambda(x^{2n})$ and $r(n,x) = \sum_{k=0}^{n} a(n-k) x^k$.

Then

$$\frac{d(n,x)}{d(n,0)} = \frac{(-1)^n}{t_{2n-1}t_{2n-3}\cdots t_1}\sum_{j=0}^{n}(-1)^k v(2n+1,k) r(n-k,x). \tag{8}$$

If $A(n)$ is defined by $A(2n) = a(n)$ and $A(2n+1) = 0$, and if $R(n,x) = \sum_{k=0}^{n} A(n-k) x^k$, then $R(2n,x) = r(n,x^2)$ and therefore $H_{2n}(0,0) = \sum_{j=0}^{n} u(j,0) R(2j,x)$ and

$H_{2n+1}(0,0) = \sum_{j=0}^{n} u(j,0) R(2j,x)$ coincide with $h_n(0,0) = \sum_{j=0}^{n} u(j,0) r(j,x^2)$.

*Therefore the corresponding quotients of Hankel determinants satisfy*

$$\frac{D(2n,x)}{D(2n,0)} = \frac{D(2n+1,x)}{D(2n+1,0)} = \frac{d(n,x^2)}{d(n,0)}. \tag{9}$$

We can now prove

**Theorem 1**

Let $t = (t_n)_{n\geq 0}$ with $t_n \neq 0$ for all $n$ and let $T = (T_n)_{n\geq 0} = (t_{n+1})_{n\geq 0}$. Let $p_n(x,t)$ be the orthogonal polynomials associated with the sequence $t$ and $p_n(x,T)$ the orthogonal polynomials associated with $T$.
Then the quotients of the Hankel determinants of the moments corresponding to $t$ are given by

$$\frac{d(n,x^2)}{d(n,0)} = \frac{(-1)^n}{t_1 t_3 \cdots t_{2n-1}} p_{2n}(x,T). \tag{10}$$

**Proof**
Since $p_n(x,t) = x p_{n-1}(x,t) - t_{n-2} p_{n-2}(x,t)$ we get

$$v(n,k) = v(n-1,k) + t_{n-2} v(n-2, k-1). \tag{11}$$

We assert that

$$\sum_{k=0}^{n}(-1)^k v(2n+1,k)r(n-k,x^2) = p_{2n}(x,T) \tag{12}$$

and

$$\sum_{k=0}^{n+1}(-1)^k v(2n+2,k)r(n+1-k,x^2) = xp_{2n+1}(x,T). \tag{13}$$

This obviously is true for $n=0$. Therefore by induction

$$\sum_{k=0}^{n}(-1)^k v(2n+1,k)r(n-k,x^2) = \sum_{k=0}^{n}(-1)^k v(2n,k)r(n-k,x^2)$$

$$+t_{2n-1}\sum_{k=0}^{n}(-1)^k v(2n-1,k-1)r(n-k,x^2) = xp_{2n-1}(x,T)$$

$$-t_{2n-1}\sum_{k=0}^{n}(-1)^k v(2n-1,k)r(n-k-1,x^2) = xp_{2n-1}(x,T) - t_{2n-1}p_{2n-2}(x,T) = p_{2n}(x,T)$$

and

$$\sum_{k=0}^{n}(-1)^k v(2n+2,k)r(n+1-k,x^2)$$

$$= \sum_{k=0}^{n}(-1)^k v(2n+1,k)r(n+1-k,x^2) + t_{2n}\sum_{k=0}^{n}(-1)^k v(2n,k-1)r(n+1-k,x^2)$$

$$= \sum_{k=0}^{n}(-1)^k v(2n+1,k)a(n+1-k) + x^2\sum_{k=0}^{n}(-1)^k v(2n+1,k)r(n-k,x^2)$$

$$-t_{2n}\sum_{k=0}^{n}(-1)^k v(2n,k)r(n-k,x^2) = x^2 p_{2n}(x,T) - t_{2n}xp_{2n-1}(x,T) = xp_{2n+1}(x,T).$$

Here we used the fact that $\sum_{k=0}^{n}(-1)^k v(2n+1,k)a(n+1-k) = 0$ by (7).

**Remarks**

The polynomials $P_n(x,t) = p_{2n}(\sqrt{x},t)$ are orthogonal with respect to the linear functional $L$ defined by $L(P_n) = [n=0]$ and satisfy $L(x^n) = a(n)$.
They satisfy the recurrence

$$P_n(x,t) = (x - S_{n-1})P_{n-1}(x,t) - U_{n-2}P_{n-2}(x,t) \tag{14}$$

with

$$\begin{aligned} S_0 &= t_0 \\ S_n &= t_{2n-1} + t_{2n} \quad \text{for } n>0 \\ U_n &= t_{2n}t_{2n+1} \end{aligned} \tag{15}$$

Theorem 1 gives only the ratio of determinants but not $d(n,0)$. But in our context $d(n,0)$ is implicitly known because $d(n,0) = U_0^n U_1^{n-1} \cdots U_{n-1}$ (cf. e.g. [5]).

Now we consider some interesting

**Examples**

1) Define the bivariate Fibonacci polynomials $F_n(x,s)$ by $F_n(x,s) = xF_{n-1}(x,s) + sF_{n-2}(x,s)$ with initial values $F_0(x,s) = 0$ and $F_1(x,s) = 1$ and consider the polynomials

$$p_n(x) = F_{n+1}(x,-1) = \sum_{k=0}^{\lfloor \frac{n}{2} \rfloor} (-1)^k \binom{n-k}{k} x^{n-2k}. \tag{16}$$

In this case $t_n = 1$ for all $n$. Therefore $T = t$.

It is well known that the moments are the Catalan numbers

$$\Lambda(x^{2n}) = C_n = \frac{1}{n+1}\binom{2n}{n}. \tag{17}$$

Thus we get

**Corollary 1 (Ömer Eğecioğlu [10])**

For $a(n) = C_n = \dfrac{1}{n+1}\binom{2n}{n}$ the Hankel determinants are

$$\frac{d(n,x^2)}{d(n,0)} = (-1)^n F_{2n+1}(x,-1) = \sum_{k=0}^{n}(-1)^{n-k}\binom{2n-k}{k}x^{2n-2k} = \sum_{j=0}^{n}(-1)^j \binom{n+j}{n-j}x^{2j}. \tag{18}$$

2) Choose $t_0 = 2$ and $t_n = 1$ for $n > 0$. Then

$$p_n(x) = L_n(x,-1) = \sum_{k=0}^{\lfloor \frac{n}{2} \rfloor} (-1)^k \frac{n}{n-k}\binom{n-k}{k} x^{n-2k} \quad \text{for } n > 0. \tag{19}$$

Here $L_n(x,s)$ are bivariate Lucas polynomials defined by $L_n(x,s) = xL_{n-1}(x,s) + sL_{n-2}(x,s)$ with initial values $L_0(x,s) = 2$ and $L_1(x,s) = x$. Note that $p_0(x) = 1 \neq L_0(x,-1) = 2$.

In this case

$$\Lambda(x^{2n}) = \binom{2n}{n} \tag{20}$$

and $T_n = 1$. Therefore we get

**Corollary 2 (Ömer Eğecioğlu [10])**

For $a(n) = \binom{2n}{n}$ the quotients of the Hankel determinants are also given by

$$\frac{d(n,x^2)}{d(n,0)} = (-1)^n F_{2n+1}(x,-1) = \sum_{k=0}^{n}(-1)^{n-k}\binom{2n-k}{k}x^{2n-2k} = \sum_{j=0}^{n}(-1)^j \binom{n+j}{n-j}x^{2j}. \quad (21)$$

3) Define the bivariate (Carlitz -) $q$ – Fibonacci polynomials $F_n(x,s,q)$ by

$$F_n(x,s,q) = xF_{n-1}(x,s,q) + q^{n-3}sF_{n-2}(x,s,q) \quad (22)$$

with initial values $F_0(x,s,q) = 0$ and $F_1(x,s,q) = 1$.

They satisfy $F_n(x,s,q) = \sum_{k=0}^{\lfloor\frac{n-1}{2}\rfloor} s^k q^{k^2-k} \begin{bmatrix} n-1-k \\ k \end{bmatrix} x^{n-1-2k}$.

Let now $t_n = q^n$. Then

$$p_n(x) = F_{n+1}(x,-1,q) = \sum_{k=0}^{\lfloor\frac{n}{2}\rfloor} (-1)^k q^{k^2-k} \begin{bmatrix} n-k \\ k \end{bmatrix} x^{n-2k}. \quad (23)$$

In this case the moments

$$\Lambda(x^{2n}) = C_n(q) \quad (24)$$

are the $q$ – Catalan numbers $C_n(q)$ of Carlitz whose generating function $f(z) = \sum_{n\geq 0} C_n(q)z^n$
satisfies $f(z) = 1 + zf(z)f(qz)$. (See e.g. [4]).

Since $t_{2n-1}t_{2n-3}\cdots t_1 = q^{1+3+\cdots+2n-1} = q^{n^2}$ and $t_{n+1} = q^{n+1}$ we have $p_n(x,T) = F_{n+1}(x,-q,q)$.

This implies

**Corollary 3**
For $a(n) = C_n(q)$ the quotients of the Hankel determinants are

$$\frac{d(n,x^2)}{d(n,0)} = \frac{(-1)^n}{q^{n^2}} F_{2n+1}(x,-q,q) = \sum_{k=0}^{n}(-1)^{n-k}q^{k^2-n^2}\begin{bmatrix} 2n-k \\ k \end{bmatrix} x^{2n-2k}. \quad (25)$$

**Remark**
Christian Krattenthaler [12] (unpublished) has previously proved Corollary 3 with another method.

4) Define the $(q,b)$ – Fibonacci polynomials $F_n(x,b,s,q)$ by the recursion

$$F_n(x,b,s,q) = xF_{n-1}(x,b,s,q) + \frac{q^{n-3}s}{(1-q^{n-2}b)(1-q^{n-1}b)} F_{n-2}(x,b,s,q) \tag{26}$$

with initial values $F_0(x,b,s,q) = 0$ and $F_1(x,b,s,q) = 1$.
These are variants of the Al Salam and Ismail polynomials ([1]).

Then (cf. e.g. [7])

$$F_n(x,b,s,q) = \sum_{k=0}^{\lfloor \frac{n-1}{2} \rfloor} s^k q^{k^2-k} \begin{bmatrix} n-1-k \\ k \end{bmatrix} \frac{x^{n-1-2k}}{\prod_{j=1}^{k}(1-q^j b) \prod_{j=n-k}^{n-1}(1-q^j b)}.$$

Let $t_n = \dfrac{q^n}{(1+q^{n+1})(1+q^{n+2})}$.

The corresponding orthogonal polynomials are

$$p_n(x,t) = F_{n+1}(x,-1,-1,q) = \sum_{k=0}^{\lfloor \frac{n}{2} \rfloor} (-1)^k q^{k^2-k} \begin{bmatrix} n-k \\ k \end{bmatrix} \frac{x^{n-2k}}{\prod_{j=1}^{k}(1+q^j) \prod_{j=n+1-k}^{n}(1+q^j)}. \tag{27}$$

Since $t_{n+1} = \dfrac{q^{n+1}}{(1+q^{n+2})(1+q^{n+3})}$ we get $p_n(x,T) = F_{n+1}(x,-q,-q,q)$.

The moments are the $q$ – Catalan numbers of George Andrews ([2])

$$\Lambda(x^{2n}) = \frac{1}{[n+1]} \begin{bmatrix} 2n \\ n \end{bmatrix} \frac{1+q}{1+q^{n+1}} \frac{1}{\prod_{j=1}^{n}(1+q^j)^2}. \tag{28}$$

A proof can be found in [6].

This implies

**Corollary 4**

*For the (Andrews-) $q$ – Catalan numbers $a(n) = \dfrac{1}{[n+1]} \begin{bmatrix} 2n \\ n \end{bmatrix} \dfrac{1+q}{1+q^{n+1}} \dfrac{1}{\prod_{j=1}^{n}(1+q^j)^2}$*

*we get*

$$\frac{d(n,x^2)}{d(n,0)} = \frac{(-1)^n}{q^{n^2}} \prod_{j=2}^{2n+1}(1+q^j)F_{2n+1}(x,-q,-q,q)$$

$$= \sum_{k=0}^{n} q^{k^2-n^2}(-1)^{n-k} \prod_{j=k+2}^{2n-k+1}(1+q^j)\begin{bmatrix}2n-k\\k\end{bmatrix}x^{2n-2k}.$$

(29)

5) Consider the generalized $q$ – Lucas polynomials (cf. [7],[8])

$$L_n(x,s,q) = \sum_{k=0}^{\lfloor\frac{n}{2}\rfloor} q^{k^2-k}s^k x^{n-2k} \frac{[n]}{[n-k]}\begin{bmatrix}n-k\\k\end{bmatrix}\frac{1}{\prod_{j=1}^{k}(1+q^j)\prod_{j=n-k}^{n-1}(1+q^j)}.$$

(30)

They satisfy

$$L_n(x,s,q) = xL_{n-1}(x,s,q) - \frac{q^{n-2}s}{(1+q^{n-2})(1+q^{n-1})}L_{n-2}(x,s,q)$$

(31)

with initial values $L_0(x,s,q) = 2$ and $L_1(x,s,q) = x.$

The corresponding orthogonal polynomials $p_n(x)$ are defined by $p_0(x) = 1$ and $p_n(x) = L_n(x,-1,q)$ for $n > 0.$

The moments are (cf. [6])

$$\Lambda(x^{2n}) = \begin{bmatrix}2n\\n\end{bmatrix}\frac{1}{\prod_{j=1}^{n}(1+q^j)^2}.$$

(32)

The corresponding $t_n$ are $t_0 = \frac{1}{1+q}$ and $t_n = \frac{q^n}{(1+q^n)(1+q^{n+1})}.$

In this case $t_1 t_3 \cdots t_{2n-1} = \frac{1}{q^{n^2}}\prod_{j=1}^{2n}(1+q^j)$ and $t_{n+1} = \frac{q^{n+1}}{(1+q^{n+1})(1+q^{n+2})}.$

Therefore $p_n(x,T) = F_{n+1}(x,-1,-q,q).$

This implies

**Corollary 5**

For $a(n) = \begin{bmatrix}2n\\n\end{bmatrix}\frac{1}{\prod_{j=1}^{n}(1+q^j)^2}$

*the quotients of the Hankel determinants are*

$$\frac{d(n,x^2)}{d(n,0)} = \frac{(-1)^n}{q^{n^2}}\prod_{j=1}^{2n}(1+q^j)F_{2n+1}(x,-1,-q,q).$$

(33)

6) For $t_{2n} = q^n a$ and $t_{2n+1} = q^n b$ the orthogonal polynomials are (cf. [4])

$$p_{2n}(x,t) = \sum_{k=0}^{n}(-a)^{n-k}q^{\binom{n-k}{2}}x^{2k}\sum_{j=0}^{n-k}\begin{bmatrix}n-j\\k\end{bmatrix}\begin{bmatrix}k+j-1\\j\end{bmatrix}\left(\frac{b}{a}\right)^{j} \tag{34}$$

and

$$p_{2n+1}(x,t) = \sum_{k=0}^{n}(-a)^{n-k}q^{\binom{n-k}{2}}x^{2k+1}\sum_{j=0}^{n-k}\begin{bmatrix}n-j\\k\end{bmatrix}\begin{bmatrix}k+j\\j\end{bmatrix}\left(\frac{b}{a}\right)^{j}. \tag{35}$$

Therefore the Hankel determinants for the moment sequence satisfy

$$\frac{d(n,x^2)}{d(n,0)} = \sum_{k=0}^{n}(-1)^{k}q^{\binom{k+1}{2}-nk}x^{2k}\sum_{j=0}^{n-k}q^{j}\begin{bmatrix}n-j\\k\end{bmatrix}\begin{bmatrix}k+j-1\\j\end{bmatrix}a^{j}b^{-k-j}. \tag{36}$$

For $q=1$ some of the moments are well-known. For $(a,b)=(1,2)$ the moments are the little Schröder numbers $(1,1,3,11,45,197,\cdots)$ and for $(a,b)=(2,1)$ we get the (large) Schröder numbers $(1,2,6,22,90,394,\cdots)$. (Cf. [6] and OEIS [13] A006318 and A001003).

**Corollary 6**

Let $(a(n))$ be the sequence of little Schröder numbers and $t=(1,2,1,2,1,2,\cdots)$. Then

$$\frac{d(n,x^2)}{d(n,0)} = \frac{(-1)^n}{2^n}p_{2n}(x,T) = \sum_{k=0}^{n}(-1)^{k}x^{2k}\sum_{j=0}^{n-k}\binom{n-j}{k}\binom{k+j-1}{j}2^{-k-j}. \tag{37}$$

If $(a(n))$ is the sequence of large Schröder numbers then

$$\frac{d(n,x^2)}{d(n,0)} = (-1)^n p_{2n}(x,t) = \sum_{k=0}^{n}(-1)^{k}x^{2k}\sum_{j=0}^{n-k}\binom{n-j}{k}\binom{k+j-1}{j}2^{j}. \tag{38}$$

**Remarks**

It is clear that for each sequence $T$ there are many sequences $t$ such that the divided Hankel determinants are $\dfrac{d(n,x^2)}{d(n,0)} = \dfrac{(-1)^n}{t_1 t_3 \cdots t_{2n-1}} p_{2n}(x,T)$.
It suffices to choose $t_0 \neq 0$ arbitrary.

It should be noted that not every sequence with non-vanishing Hankel determinants can be represented as moments of symmetric orthogonal polynomials. For example the sequence $(1,1,2,4,9,21,51,127,\cdots)$ of Motzkin numbers $M_n$ (cf. OEIS [13], A001006) satisfies $\det(M_{i+j})_{i,j=0}^n = 1$ for all $n$. The first divided Hankel determinants $\frac{d(n,x^2)}{d(n,0)}$ for the Motzkin numbers turn out to be $1,\ 1-x^2,\ 1-x^2,\ 1-x^2-2x^4+x^6,\cdots$.

Thus $\frac{d(2,x^2)}{d(2,0)}$ is a polynomial of degree 2 instead of degree 4. Thus $M_n$ cannot be of the form $\Lambda(x^{2n})$ corresponding to orthogonal polynomials of the form (2).

So what can be said about the ratios of Hankel determinants of the Motzkin numbers? To answer this question we need the following fact (cf. e.g. [4]). Let $f(z,u) = \sum_{n\geq 0} M_n(u)z^n$ satisfy $f(z,u) = 1 + uzf(z,u) + z^2 f(z,u)^2$ and define the linear functional $F_u$ by $F_u(x^n) = M_n(u)$. Then the corresponding orthogonal polynomials satisfy

$$P_n(x,t) = (x-u)P_{n-1}(x,t) - P_{n-2}(x,t). \tag{39}$$

Thus

$$P_n(x,t) = F_{n+1}(x-u,-1). \tag{40}$$

For $u=1$ we get the Motzkin numbers $M_n(1) = M_n$. In this case there are no $t_n$ satisfying (15). But for general $u$ it is easily seen that $t_{2n} = \frac{F_{n+2}(u,-1)}{F_{n+1}(u,-1)}$ and $t_{2n+1} = \frac{1}{t_{2n}}$.

Therefore (10) implies

$$\frac{d(n,x^2)}{d(n,0)} = (-1)^n F_{n+1}(u,-1) p_{2n}(x,T). \tag{41}$$

Now observe that the polynomials $p_{2n+1}(x,T)$ satisfy both
$p_{2n+1}(x,T) = xp_{2n}(x,T) - t_{2n} p_{2n-1}(x,T)$
and
$p_{2n+3}(x,T) = (x^2 - u) p_{2n+1}(x,T) - p_{2n-1}(x,T)$.

This implies that

$$p_{2n}(x,T) = F_{n+1}(x^2-u,-1) + \frac{F_{n+2}(u,-1)}{F_{n+1}(u,-1)} F_n(x^2-u,-1).$$

Therefore (41) gives

$$\frac{d(n,x)}{d(n,0)} = (-1)^n \left( F_{n+1}(u,-1) F_{n+1}(x-u,-1) + F_{n+2}(u,-1) F_n(x-u,-1) \right). \tag{42}$$

By continuity this relation also holds for $u = 1$. Observing that the sequence $(F_n(1,-1))_{n \geq 0} = (0,1,1,0,-1,-1,\cdots)$ is periodic with period 6 we see that in this case

$$d(3n, x) = F_{3n}(x-1,-1) + F_{3n+1}(x-1,-1) \tag{43}$$

and

$$d(3n+1, x) = d(3n+2, x) = -F_{3n+2}(x-1,-1). \tag{44}$$

**2. The case $r(n,x) = a(n)x - a(n+1)$.**

Generalizing the results of Aleksandar Cvetković, Predrag Rajković and Miloš Ivković [9] I have shown in [5] a general result which in the present terminology can be stated in the following way:

**Theorem 2**

*Let*

$$r(n, x) = a(n)x - a(n+1). \tag{45}$$

*Then*

$$\frac{\det\left(r(i+j, x^2)\right)_{i,j=0}^n}{\det\left(a(i+j)\right)_{i,j=0}^n} = p_{2n+2}(x,t) \tag{46}$$

*and*

$$\frac{\det\left(r(i+j+1, x^2)\right)_{i,j=0}^n}{\det\left(a(i+j+1)\right)_{i,j=0}^n} = \frac{p_{2n+3}(x,t)}{x}. \tag{47}$$

**Proof**

Here I give a simpler proof. With the notation introduced in (3) we get

$$x^{2m} p_{2n+2}(x) = \sum_{k=0}^{n+1} (-1)^k v(2n+2,k) x^{2n+2+2m-2k}.$$

Applying the linear functional $\Lambda$ to this identity we get

$$a(n+1+m) = \sum_{k=1}^{n+1} (-1)^{k+1} v(2n+2,k) a(n+1+m-k) = \sum_{j=0}^{n} (-1)^{n-j} v(2n+2, n+1-j) a(m+j)$$

for $0 \leq m \leq n$.

This implies that

$$\left(r(i+j, x^2)\right)_{i,j=0}^n = \left(b(i,j)\right)_{i,j=0}^n \left(a(i+j)\right)_{i,j=0}^n$$

with

$$(b(i,j))_{i,j=0}^{n} = \begin{pmatrix} x^2 & -1 & 0 & \cdots & 0 \\ 0 & x^2 & -1 & \cdots & 0 \\ \vdots & \vdots & \vdots & \ddots & \vdots \\ 0 & 0 & 0 & \ddots & -1 \\ (-1)^n v(2n+2,n+1) & (-1)^{n-1} v(2n+2,n) & (-1)^n v(2n+2,n-1) & \cdots & x^2 - v(2n+2,1) \end{pmatrix}.$$

Therefore

$$\det(b(i,j))_{i,j=0}^{n} = \sum_{k=0}^{n+1} (-1)^k v(2n+2,k) x^{2n+2-2k} = p_{2n+2}(x).$$

Identity (47) follows analogously from

$$x^{2m+2} \frac{p_{2n+3}(x)}{x} = \sum_{k=0}^{n+1} (-1)^k v(2n+3,k) x^{2n+4+2m-2k},$$

which for $0 \leq m \leq n$ implies

$$a(n+2+m) = \sum_{j=0}^{n} (-1)^{n-j} v(2n+3, n+1-j) a(m+1+j).$$

Note that (46) for $x = 0$ implies that $\det(a(i+j+1))_{i,j=0}^{n} \neq 0$ since
$p_{2n+2}(0,t) = (-1)^{n+1} t_0 t_2 \cdots t_{2n} \neq 0.$

**Examples**

1) For $a(n) = C_n$ we have $p_n(x,t) = F_{n+1}(x,-1)$.

This implies the results obtained by Aleksandar Cvetković, Predrag Rajković and Miloš Ivković [9].

Observe that $F_n(\sqrt{-1}x, -1) = (\sqrt{-1})^{n-1} F_n(x,1).$
Therefore (46) reduces to

$$\frac{\det(r(i+j, -x^2))_{i,j=0}^{n}}{\det(a(i+j))_{i,j=0}^{n}} = F_{2n+3}(\sqrt{-1}x, -1) = (-1)^{n+1} F_{2n+3}(x,1).$$

This implies for $x = 0$ that $\det(C_{i+j+1})_{i,j=0}^{n} = 1$

and for $x = 1$ that

$$\det(C_{i+j} + C_{i+j+1})_{i,j=0}^{n} = F_{2n+3}(1,1) = F_{2n+3}, \qquad (48)$$

where $F_n$ denotes a Fibonacci number.

In the same way (47) gives

$$\det\left(C_{i+j+1} + C_{i+j+2}\right)_{i,j=0}^{n} = F_{2n+4}. \tag{49}$$

2) As another simple example consider the Hermite polynomials $H_n(x)$ defined by
$H_n(x) = xH_{n-1}(x) - (n-1)H_{n-2}(x)$.
The moments are $\Lambda(x^{2n}) = (2n-1)!!$ (cf. e.g. [14]).
Therefore we get

$$\frac{\det\left((2i+2j-1)!!(x^2 - 2i - 2j - 1)\right)_{i,j=0}^{n}}{\det\left((2i+2j-1)!!\right)_{i,j=0}^{n}} = H_{2n+2}(x) \tag{50}$$

and

$$\frac{\det\left((2i+2j+1)!!(x^2 - 2i - 2j - 3)\right)_{i,j=0}^{n}}{\det\left((2i+2j+1)!!\right)_{i,j=0}^{n}} = \frac{H_{2n+3}(x)}{x}. \tag{51}$$

3) More interesting are the numbers $M_n(u)$. The determinants $\det\left(M_{i+j}(u) + M_{i+j+1}(u)\right)$ and $\det\left(M_{i+j+1}(u) + M_{i+j+2}(u)\right)$ have been computed in [3].
These results are also simple consequences of Theorem 2.

By (40) we know that $p_{2n}(\sqrt{x}, t) = P_n(x, t) = F_{n+1}(x - u, -1)$.
Therefore

$$\frac{\det\left(xM_{i+j}(u) - M_{i+j+1}(u)\right)_{i,j=0}^{n}}{\det\left(M_{i+j}(u)\right)_{i,j=0}^{n}} = F_{n+2}(x - u, -1). \tag{52}$$

For $x = 0$ this gives

$$\det\left(M_{i+j+1}(u)\right)_{i,j=0}^{n} = (-1)^{n+1} F_{n+2}(-u, -1) = F_{n+2}(u, -1). \tag{53}$$

For $x = -1$ we get from (52)

$$\det\left(M_{i+j}(u) + M_{i+j+1}(u)\right)_{i,j=0}^{n} = F_{n+2}(u+1, -1). \tag{54}$$

Note that for $u = 1$ the sequence $\left(F_n(1, -1)\right)_{n \geq 0} = (0, 1, 1, 0, -1, -1, \cdots)$ is periodic with period 6 and $F_n(2, -1) = n$.

For the second formula observe that $p_{2n+1}(x, t) = xp_{2n}(x, t) - t_{2n-1}p_{2n-1}(x, t)$ and therefore
$F_{n+2}(u, -1)p_{2n+3}(x, t) = F_{n+2}(u, -1)p_{2n+2}(x, t) - F_{n+1}(u, -1)p_{2n+1}(x, t)$ which gives by iteration

$$F_{n+2}(u,-1)p_{2n+3}(x,t) = \sum_{j=0}^{n+1}(-1)^j F_{n+2-j}(u,-1)p_{2n+2-2j}(x,t) = \sum_{k=0}^{n+1}(-1)^{n+1-k}F_{k+1}(u,-1)F_{k+1}(x^2-u,-1).$$

Therefore we get

$$\frac{\det\left(xM_{i+j+1}(u) - M_{i+j+2}(u)\right)_{i,j=0}^{n}}{\det\left(M_{i+j+1}(u)\right)_{i,j=0}^{n}} = \frac{(-1)^{n+1}}{F_{n+2}(u,-1)}\sum_{k=0}^{n+1}(-1)^k F_{k+1}(u,-1)F_{k+1}(x-u,-1). \quad (55)$$

For $u = 1$ this does not make sense for all $n$ since $F_{3n}(1,-1) = 0$. But since for general $u$ $\det\left(M_{i+j+1}(u)\right)_{i,j=0}^{n} = F_{n+2}(u,-1)$ we get by multiplying (55) with $F_{n+2}(u,-1)$

$$\det\left(xM_{i+j+1}(u) - M_{i+j+2}(u)\right)_{i,j=0}^{n} = \sum_{k=0}^{n+1}(-1)^{n+1-k}F_{k+1}(u,-1)F_{k+1}(x-u,-1). \quad (56)$$

Since both sides are polynomials in $u$ this also holds for $u = 1$.

4) Finally I state the corresponding results for the little Schröder numbers $s(n)$ which follow from (34) and (35). These results have also been obtained in [11] with another method.

$$\frac{\det\left(xs(i+j) - s(i+j+1)\right)_{i,j=0}^{n}}{\det\left(s(i+j)\right)_{i,j=0}^{n}} = \sum_{k=0}^{n+1}(-1)^{n+1-k}x^k \sum_{j=0}^{n+1-k} 2^j \binom{n+1-j}{k}\binom{k+j-1}{j}. \quad (57)$$

The first terms are $x - 1, x^2 - 4x + 1, x^3 - 7x^2 + 11x - 1, \cdots$.

$$\frac{\det\left(xs(i+j+1) - s(i+j+2)\right)_{i,j=0}^{n}}{\det\left(s(i+j+1)\right)_{i,j=0}^{n}} = \sum_{k=0}^{n+1}(-1)^{n+1-k}x^k \sum_{j=0}^{n+1-k} 2^j \binom{n+1-j}{k}\binom{k+j}{j}. \quad (58)$$

The first terms are $x - 3, x^2 - 6x + 7, x^3 - 9x^2 + 23x - 15, \cdots$.